\documentclass[11pt,a4paper, twoside]{article}
\usepackage{amssymb}
\usepackage{mathrsfs}
\usepackage{amsfonts}
\usepackage{color}
\usepackage{mathbbol}
\usepackage{amsmath,amsthm,environ,bbm}

\allowdisplaybreaks

 \NewEnviron{ews}{%
\begin{equation}\begin{split}
  \BODY
\end{split}\end{equation}
}

\NewEnviron{ews*}{%
\begin{equation*}\begin{split}
  \BODY
\end{split}\end{equation*}
}

\def\beg{\begin}
\def\bequ{\begin{equation}}
\def\enqu{\end{equation}}
\def\bes{\begin{split}}
\def\ens{\end{split}}
\def\bews{\begin{ews}}
\def\beqn{\begin{eqnarray}}
\def\enqn{\end{eqnarray}}
\def\beq*{\begin{equation*}}
\def\enq*{\end{equation*}}
\def\bqn*{\begin{eqnarray*}}
\def\eqn*{\end{eqnarray*}}
\def\bary{\begin{array}}
\def\eary{\end{array}}
\def\bpma{\begin{pmatrix}}
\def\epma{\end{pmatrix}}
\def\bvma{\begin{Vmatrix}}
\def\evma{\end{Vmatrix}}

 \newtheorem{thm}{Theorem}[section]

 \numberwithin{equation}{section}

\def\al{\alpha}
\def\be{\beta}
\def\ga{\gamma}
\def\de{\delta}

\def\th{\theta}

\def\si{\sigma}

\def\De{\Delta}

\def\R{\mathbb R}

\def\Z{\mathbb Z}

\def\d{\mathrm{d}}

\def\ff{\frac}
\def\ra{\rightarrow}
\def\nn{\nabla}

\def\<{\langle}
\def\>{\rangle}

\def\1{\mathbbm{1}}

\makeatletter
\renewcommand\@biblabel[1]{${#1}.$}
\makeatother

\def\vd{\mathrm{d}}

\def\bthm{\begin{theorem}}
\def\ethm{\end{theorem}}
\def\blem{\begin{lemma}}
\def\elem{\end{lemma}}
\def\brem{\begin{remark}}
\def\erem{\end{remark}}
\def\bexm{\begin{example}}
\def\eexm{\end{example}}

\def\ecor{\end{corollary}}
\def\r{\right}
\def\l{\left}

\def\e{\operatorname{e}}

\def\be{\begin{equation}}

\def\zm{\noindent{\bf  Proof.\ }}
\def\endzm{\quad $\Box$}

\def\Rd{\mathbb{R}^d}

\newtheorem{theorem}{Theorem}[section]
\newtheorem{corollary}[theorem]{Corollary}
\newtheorem{lemma}[theorem]{Lemma}

\newtheorem{example}[theorem]{Example}
\theoremstyle{definition}

\newtheorem{remark}[theorem]{Remark}

\numberwithin{equation}{section}

\def\ff{\frac}
\def\ra{\rightarrow}

\begin{document}

\arraycolsep=1pt

\title{\bf Weak Poincar\'{e} Inequalities for Convolution Probabilities Measures
\footnotetext{\hspace{-0.35cm} 2010 {\it Mathematics Subject
Classification}. {Primary 60J75; Secondary   47G20; 60G52.}
\endgraf{\it Key words and phrases}. Weak Poincar\'{e} inequality, Lyapunov condition, convolution.}}
\author{Li-Juan Cheng\,  and  Shao-Qin Zhang \footnote{Corresponding author}}
\date{}
\maketitle

\vspace{-0.5cm}

\begin{center}
\begin{minipage}{13.5cm}\small
{{\bf Abstract.} By using Lyapunov conditions,  weak Poincar\'{e} inequalities are established for some probability measures on a manifold $(M,g)$. These results are further applied to the convolution of two probability measures on $\R^d$. Along with explicit results we study concrete examples.
  }
\end{minipage}
\end{center}

\baselineskip 17pt
%%%%%%%%%%%%%%%%%%%%%%%%%%%%%%%%%%%%%%%%%%%%%%%%%%%%%%%%%%%%%%%%%%%%%%%%%%%%%%%%%%%%

\section{Introduction}
During the last decades, a lot of attention has been devoted to the study of ergodic theory
for Markov processes. Specifically a lot of effort has been made on the stability speed for the corresponding Markov processes (see e.g. \cite{BCG,MT93a,MT93b,MT93c,Wbook1}). From this former
work,  functional inequalities of Dirichlet forms play important roles in characterizing
the convergence speed of ergodic Markov processes. For instance, Poincar\'{e} inequalities imply
the exponential ergodic speed of Markov processes; super Poincar\'{e} inequalities imply the strong
ergodicity of the corresponding processes; weak Poincar\'{e} inequalities are used to characterize
the non-exponential convergence rate for semigroup (see \cite{Wbook1} for details).

However,  to establish a functional inequality, we always need the coefficients of the  generator to  satisfy some regularity conditions. To deal with generators with less regular or less explicit coefficients,
an efficient way is to regard the measures as perturbations from better ones, which satisfy the
underlying functional inequalities. The convolution probability measure, in the sense of an independent
sum of random variables, can be regarded as a kind of perturbation; see e.g. \cite{Cha8,Zim}
and references therein. Moveover, the study of functional inequalities for convolution probability
measures is helpful in describing some behaviors of random variables under independent
perturbations, see e.g. \cite[Section 3]{Zim} for an application to the study of random matrices.

Recently, F.-Y. Wang and J. Wang \cite{WW} gave some sufficient conditions for  log-Sobolev/ Poincar\'{e}/ super Poincar\'{e} inequalities for convolution probability measures.   The
present article is thus a continuation of \cite{WW} to  study  weak  Poincar\'{e} inequalities for the convolution probability measures.

Before moving on, let us briefly review some background about the weak Poincar\'{e} inequality. The weak Poincar\'{e} inequality was first introduced in \cite{RW} to characterize the non-exponential convergence rate
of Markov processes and  the concentration of measure phenomenon for sub-exponential laws  (see \cite{BCR}). Let $(M,g)$ be a $d$-dimensional complete connected Riemannian manifold and $\d x$ be the volume measure. For a probability measure $\mu(\vd x):=\e^{-V(x)}\,\vd x$ with some locally bounded function
$V$ on $M$, we say that $\mu$ satisfies the \emph{weak Poincar\'{e} inequality} if
 \begin{equation}\label{weak-poin}
 \|f\|^2\leq \alpha(r)\mu(|\nabla f|^2)+r{\rm Osc}^2(f), \ \ r>0, f\in C_b^2(M)
 \end{equation}
 holds for some decreasing function $\alpha:[0,\infty)\rightarrow (0,\infty)$,  where $\|\cdot\|$ denotes the $L^2(\mu)$-norm and ${\rm Osc}(f):=\sup_{x,y\in M}|f(x)-f(y)|$. Indeed, the function $\alpha$ can be estimated by using the growth of
 $|V|$ (see \cite{RW,Wbook1}). However, in general, the resulting estimate of the rate function is less sharp. Therefore, in Section 2, we will revisit this problem on Riemannian manifolds by using some Lyapunov conditions.

As an application of the results in Section 2, we consider the weak Poincar\'{e} inequality for convolution probability measures on $\R^d$.
Let $\mu$ and $\nu$ be two probability measures on $\mathbb{R}^d$. The perturbation of $\mu$ by the probability measure $\nu$ is given by their convolution
 $$(\nu*\mu)(A):=\iint_{\mathbb{R}^d\times \mathbb{R}^d}\mathbb{1}_A(x+y)\mu(\vd x)\nu(\vd y),$$
 where $A\in \mathcal{B}(\mathbb{R}^d)$. In particular, let $\mu(\vd x)=\e^{-V(x)}\,\vd x$ be a probability measure on $\mathbb{R}^d$ such that $V\in C^1(\mathbb{R}^d)$ and $\nu$ be a probability measure
on $\mathbb{R}^d$ such that
\bequ\label{equ-e0}
p_{\nu}(\cdot):=\int\e^{-V(\cdot-z)}\,\nu(\vd z)\in C^1(\R^d).
\enqu
Then
\begin{equation}\label{in-e1}
(\mu\ast \nu)(\vd x)=p_{\nu}(x)\,\vd x=\e^{-V_{\nu}(x)}\,\vd x,
\end{equation}
where $V_{\nu}(x):=-\log p_{\nu}(x)$.
Let $L_{\nu}=\Delta-\nabla V_{\nu}$, which is the generator associated with some independent sum of two Markov diffusion processes with invariant
measures $\mu$ and $\nu$, respectively.  This article aims to prove that the measure $\mu*\nu$ satisfies \eqref{weak-poin} for some explicit function $\alpha$, which characterizes the explicit $L^2$-ergodic speed of some diffusion generated by $L_{\nu}$. Actually, the existence of weak Poincar\'{e} inequalities for $\mu*\nu$ holds automatically due to the positivity of  the density  $\e^{-V_{\nu}(x)}$ (see \cite{RW}). So the main topic of this article is to find an explicit function $\alpha$ in the weak Poincar\'{e} inequality.

Our method is based on the use of Lyapunov type conditions. These conditions are well known
to furnish some results on the long time behavior of the laws of Markov processes (see e.g.
\cite{BCG,CGWW,MT93a,MT93b,MT93c,WW} and references therein). In the recent work \cite{WW}, the authors partly use
Lyapunov conditions to study ordinary or super Poincar\'{e} inequality for convolution probability
measures. As announced, the present paper is thus a complement of \cite{WW} for the study of the weak
Poincar\'{e} inequality. The main idea of the use of a Lyapunov function is similar to \cite{WW} and in
the present work, however we have to face some technical difficulties when choosing
suitable Lyapunov functions and handling the ``local term" in the proof of Theorem \ref{th-1} below.
It is worthy to mentioning that a new and reasonable Lyapunov function, constructed for establishing weak Poincar\'{e} inequalities,  can also be applied to improving some results obtained in \cite{WW} for the super Poincar\'{e} inequality. In addition we will use a
comparison method  to simplify the assumptions in general results, and then give some
{concrete} examples as applications.

The parts of the paper are organized as follows. In the following section, we study the weak Poincar\'{e}
 inequality by  Lyapunov conditions and the comparison theorem for some probability measures on Riemannian manifolds. In Section 3,  we apply results in Section 2 to convolution probability measures on $\R^d$.  Some explicit examples are studied in Section 4.

\section{Weak Poincar\'{e} inequality on manifolds via Lyapunov conditions}
We organize this section by first introducing main results and then giving proofs.
\subsection{Main results}
Let $(M,g)$ be a $d$-dimensional complete connected Riemannian manifold.
Let $\nabla$ and $\Delta$ be the Levi-Civita connection and the Laplacian associated with $g$, respectively. Consider the elliptic operator $L=\Delta-\nabla V$ for $V\in C^1(M)$ such that
$\mu(\vd x):=\e^{-V(x)}\,\vd x$ is a probability measure, where $\,\vd x$ is the Riemannian volume measure.

Given $o\in M$. For any $x\in M$, let $\rho_o(x)$ be the Riemannian distance on $M$ between $x$ and $o$ and ${\rm Cut}_o$ be the set of cut-locus points of $o$ {which is closed and has volume zero}. Define $\varphi(s)$ to be the continuous version of $$\inf_{\rho_o(x)=s,\ x\notin {\rm Cut}_o}(\l<\nabla V(x), \nabla \rho_o(x)\r>-\Delta \rho_o(x))$$ for $s\geq 0$.
 We now introduce the main results about  weak Poincar\'{e} inequalities for $\mu$ via Lyapunov conditions.
\begin{theorem}\label{th-1}
Let  $\mu(\vd x)=\e^{-V(x)}\,\vd x$ be a probability measure on $M$ for some $V\in C^{1}(M)$.
\begin{itemize}
  \item [{\rm (a)}]
Assume that for some constant $R_0$ and any $\sigma \in (0,1)$, one has
      $$\theta(r):=\ff {(1-\sigma)\varphi(r)\exp\big[\sigma\int_{R_0}^{r}\varphi(u)\,\vd u\big]}{ \int_{R_0}^{r}\exp\big[\sigma\int_{R_0}^s\varphi(u)\,\vd u\big]\,\vd s+1}>0,\ \  r\geq R_0.$$
     Let \begin{align}\label{eq-4}\phi(x)=\theta(\rho_o(x)\vee R_0),\quad x\in M.\end{align}
     Then $\mu$ satisfies the weak Poincar\'{e} inequality with $\alpha(r):=cF_{\phi}^{-1}(r)$ for some positive constant $c$, where $F_{\phi}(r):=\mu(\phi\leq \frac{1}{r})$ and $F_{\phi}^{-1}(r)=\inf\{s:~F_{\phi}(s)\leq r\}$.
  \item [{\rm (b)}] Let $V\in C^2(M)$ such that for some positive constants $R_0$ and $\delta\in (0,1)$, there exists some positive function $\phi$
  on $M$ such that
   $$\phi(x)=(1-\delta)(\delta |\nabla V|^2(x)-\Delta V(x))>0,\ \ \rho_o(x)\geq R_0.$$
   Then  $\mu$ satisfies the weak Poincar\'{e} inequality with $\alpha(r):=cF_{\phi}^{-1}(r)$ for some positive constant $c$.
\end{itemize}
\end{theorem}
\begin{remark}\label{rem-add-1}
\begin{itemize}
  \item [(i)]In Theorem \ref{th-1} (a), it is easy to see that a different $\sigma\in (0,1)$ does not affect the sign of $\theta$. But suitable choosing of $\sigma$ seems to get the best
  $\alpha$ in the  weak Poincar\'{e} inequality; see the proof of Example \ref{add-ex-2} for more explanations.
  \item [(ii)]In the proof of this theorem, we use two ways to construct Lyapunov functions, the first is new and the second is due to \cite{WW}. Our new Lyapunov function  can improve the result in \cite[Theroem 4.1(a)]{WW} for the super Poincar\'{e} inequality of convolution probability measures on $\R^d$, see Remark \ref{Rem-3} for details.
\end{itemize}
\end{remark}
We now assume that $(M,g)$ satisfies the following curvature condition:
\begin{description}
  \item [{\bf Assumption (A)}]:\  Ric$\geq -({d}-1)k$ for some constant $k$, where Ric is the Ricci curvature
tensor with respect to $g$.
\end{description}
\noindent
 Let
\begin{align*}
h_k(r)=\left\{
         \begin{array}{ll}
           \sin(\sqrt{-k r})/\sqrt{-k},\quad &\ \hbox{if}\ \  k<0; \\
           r,\quad &\ \hbox{if}\ \ k=0; \\
           \sinh(\sqrt{k}r)/\sqrt{k},\quad &\ \hbox{if}\ \ k>0.
         \end{array}
       \right.
\end{align*}
Under assumption {\bf (A)},   we can use the following comparison theorem to handle $\Delta \rho_o$ (see \cite[Section 1]{Cha8}):
    $$\Delta \rho_o\leq \frac{(d-1)h_k'(\rho_o)}{h_k(\rho_o)}$$
    outside the cut-locus. Then, the following corollary can be proved by a similar discussion as in Theorem \ref{th-1}.
\begin{corollary}\label{th-2}
Let $\mu(\vd x)=\e^{-V(\rho_o(x))}\,\vd x$ be a probability measure on $M$ for some function $V\in C^{1}(\R^+)$. Suppose that assumption {\bf (A)} holds, then we have the following two assertions.
\begin{itemize}
  \item [{\rm (a)}] Assume that for some positive constant $R_0$ and any $\sigma\in (0,1)$, one has
      $$\theta(r):=\ff {(1-\sigma)h_k(r)^{1-d}V'(r)\e^{\sigma V(r)}}{\int_{R_0}^{r}h_k(s)^{1-d}\e^{\sigma V(s)}\,\vd s+1}>0,\ \  r\geq R_0.$$
      Then $\mu$ satisfies the weak Poincar\'{e} inequality with $\alpha(r):=cF_{\phi}^{-1}(r)$ for some
      positive constant $c$, where $\phi(x):=\theta(\rho_o(x))$ for $\rho_o(x)\geq R_0$.
  \item [{\rm (b)}]  Assume that for some positive constants $R_0$ and $\delta\in (0,1)$, one has ${V}\in C^2([R_0,\infty))$ and
      \begin{align}\label{assumption-b}
   \theta(r)=(1-\delta)\l[\l(\delta -\frac{(d-1)h_k'(r)}{h_k(r)}\r){V}'(r)-{V}''(r)\r]>0,\ \ r\geq R_0.
      \end{align}
   Then  $\mu$ satisfies a weak Poincar\'{e} inequality with $\alpha(r):=cF_{\phi}^{-1}(r)$ for some positive constant $c$, where $\phi(x):=\theta(\rho_o(x))$ for $\rho_o(x)\geq R_0$.
      \end{itemize}
\end{corollary}

\subsection{Proofs}

Let $L$ be a second order elliptic operator. To prove these results above, let us first introduce the following \emph{general Lyapunov condition} with respect to $L$ (see \cite[Subsection 3.3]{DFG}).
\begin{description}
  \item [{\bf Hypothesis}\  {\bf (L)}] There exist some positive  constants $ b, r_0$,  some positive function $\phi$ on $M$ and function $W\in \mathcal{D}(L)$ with $W\geq 1$ such that
 \begin{align}\label{Lcondition}
\frac{LW}{W}\leq -\phi+b\mathbb{1}_{B_{r_0}},
 \end{align}
where $\mathcal{D}(L)$ is the weak domain of $L$ and $B_{r_0}:=\{x\in M: \rho_o(x)\leq r_0\}$ is the ball  with center $o$ and radius $r_0$.
\end{description}
\noindent Our first step is to prove that if hypothesis {\bf (L)} holds for $L=\Delta-\nabla V$, then there exists  some function $\alpha$ such that the  weak Poincar\'{e} inequality holds for $\mu$.
%\textcolor{red}{We remark here that the Lyapunov functions constructed in the following are smooth enough they all are in $\mathcal{D}(L)$.}

\begin{lemma}\label{lem1}
Let $\mu(\vd x)=\e^{-V(x)}\,\vd x$ be a probability measure on $M$. Assume the Lyapunov
condition {\bf (L)} holds for $L=\Delta-\nabla V$. Then  the following weak Poincar\'{e} inequality
$$\mu(f^2)\leq c_0F_{\phi}^{-1}(r) \mu(|\nabla f|^2)+r{\rm Osc}(f)^2$$
holds for some positive constant $c_0$ and $F_{\phi}(r):=\mu(\phi\leq \frac{1}{r})$.
\end{lemma}
\zm
The proof is given by combining
\cite[Theorem 4.6]{CGGR} with \cite[Theorem 2.18]{CGGR}. For the sake of  completeness, we include it here.
For any $r>0$ and $f\in C_b^1(M)$ with $\mu(f)=0$,  we have
\begin{align}\label{add-1e2}
\mu(f^2)=\inf_{c\in \mathbb{R}}\mu(f-c)^2&\leq\int_{\{\phi> 1/r\}}(f-f(x_0))^2\,\vd \mu+\int_{\{\phi\leq  1/r\}}(f-f(x_0))^2\,\vd \mu\nonumber\\
&\leq \int_{\{\phi> 1/r\}}(f-f(x_0))^2\,\vd \mu+\mu(\phi\leq  1/r){\rm Osc}(f)^2\nonumber\\
&\leq r\int \phi (f-f(x_0))^2\,\vd \mu+\mu(\phi \leq  1/r){\rm Osc}(f)^2\nonumber\\
&\leq -r \int\frac{LW}{W}(f-f(x_0))^2\,\vd \mu+rb\int_{B_{r_0}}(f-f(x_0))^2\,\vd \mu\nonumber\\
&\quad+\mu(\phi \leq  1/r){\rm Osc}(f)^2,
\end{align}
where  $x_0\in M$ will be specified later. Now we need to estimate the first two terms on the right hand side of the latter inequality: a {global} term and a local term. For the {global} term, by \cite[Lemma 2.12]{CGWW},  we have
\begin{align}\label{estimate-11}
-\int\frac{LW}{W}(f-f(x_0))^2\,\vd \mu\leq \int |\nabla f|^2\,\vd \mu.
\end{align}
 For the local one,
 choose $x_0\in M$ such that $f(x_0)=\frac{1}{\mu(B_{r_0})}\big(\int_{B_{r_0}}f \,\vd \mu\big)$ and define $g=f-f(x_0)$. Then we obtain
\begin{align}\label{1e1}
\int_{B_{r_0}}(f-f(x_0))^2\,\vd \mu&=\int_{B_{r_0}}g^2\,\vd \mu\leq \lambda_{r_0}^{-1}\int|\nabla g|^2\,\vd \mu +\frac{1}{\mu(B_{r_0})}\l(\int_{B_{r_0}}g \,\vd \mu\r)^2\notag\\
&=\lambda_{r_0}^{-1} \mu(|\nabla g|^2)=\lambda_{r_0}^{-1} \mu(|\nabla f|^2),
\end{align}
where by \cite[(4.3.5)]{Wbook1},
$$\lambda_{r_0}^{-1}\leq\frac{4r_0^2}{\pi^2}\exp\l\{\sup_{x,y\in B_{r_0}}(V(x)-V(y))\r\}<\infty.$$
Now, taking \eqref{estimate-11} and \eqref{1e1}    into \eqref{add-1e2}, we arrive at
\begin{align}
\mu(f^2)\leq r(b \lambda_{r_0}^{-1}+1)\int |\nabla f|^2\,\vd \mu+\mu\l(\phi\leq \frac{1}{r}\r){\rm Osc}(f)^2.
\end{align}
Let
$$F_{\phi}(r)=\mu\l(\phi \leq \frac{1}{r}\r).$$
Then
$$\lim_{r\ra+\infty}F_{\phi}(r)=0$$
due to the fact that $\phi$ is positive and $\mu$ is a probability measure on $M$. From this, we derive that
 $F_{\phi}:(0,+\infty)\rightarrow (0,1)$ is a decreasing function.
Then
$$\alpha(r):=(b\lambda_{r_0}^{-1}+1)F_{\phi}^{-1}(r)$$
is a function from $(0,+\infty)$ to $(0,+\infty)$ and the weak Poincar\'e inequality holds for such $\alpha$.
\endzm

\begin{proof}[{\bf Proof of Theorem \ref{th-1}}]
In case (a), let $0<\sigma<1$, and define the Lyapunov function by
$$W_{\sigma}(r)=\int_{R_0}^{r}\exp\bigg[\sigma\int_{R_0}^s\varphi(u)\,\vd u\bigg]\,\vd s+1,~~\mbox{for all}~~r\geq R_0.$$
By an approximation argument, we may consider $\rho_o\in C^2(M)$ for the sake of conciseness.
Then for all $\rho_o(x)\geq R_0$, we have
\begin{align*}
\frac{LW_{\sigma}(\rho_o(x))}{W_{\sigma}(\rho_o(x))}&=\ff{1}{W_{\si}(\rho_o(x))}\l[W'_{\si}(\rho_o(x))\Delta \rho_o(x)+W_{\si}''(\rho_o(x))|\nabla \rho_o(x)|^2-W_{\si}'(\rho_o(x))\l<\nabla V, \nabla \rho_o(x)\r>\r]\\
&=\ff {1}{W_{\si}(\rho_o(x))}\l\{W_{\si}'(\rho_o(x))[\Delta \rho_o(x)-\l<\nabla V,\nabla \rho_o(x)\r>]+W_{\si}''(\rho_o(x))\r\}\\
&=\ff {1}{W_{\si}(\rho_o(x))}\exp\l[\sigma\int_{R_0}^{\rho_o(x)}\varphi(u)\,\d u\r][\Delta \rho_o(x)-\l<\nabla V, \nabla \rho_o(x)\r>+\sigma \varphi(\rho_o(x))]\\
&\leq -\ff {1-\sigma}{W_{\si}(\rho_o(x))}\exp\l[\sigma\int_{R_0}^{\rho_o(x)}\varphi(u)\,\d u\r]\varphi(\rho_o(x)).
\end{align*}
Thus, there exists a constant $b>0$ such that
\begin{align*}
\frac{LW_{\sigma}(\rho_o(x))}{W_{\sigma}(\rho_o(x))}\leq& -\theta(\rho_o(x))\mathbb{1}_{\{\rho_o(x)\geq R_0\}}+b\mathbb{1}_{\{\rho_o(x)< R_0\}},
\end{align*}
which combining with Lemma \ref{lem1} implies Theorem \ref{th-1}(a).

In case (b), we consider a function $W$ in $C^2(M)$ such that $W(x)=\e^{(1-\delta)V(x)}$ for all $\rho_o(x)\geq R_0$. It is easy to see that $W(x)\geq 1$ for all $x\in M$ and
\begin{align*}
\frac{LW(x)}{W(x)}&\leq-(1-\delta)(\delta |\nabla V|^2-\Delta V)(x)\mathbb{1}_{\{\rho_o(x)\geq R_0\}}+b\mathbb{1}_{\{\rho_o(x)<R_0\}}\\
&=-\phi(x)\mathbb{1}_{\{\rho_o(x)\geq R_0\}}+b\mathbb{1}_{\{\rho_o(x)< R_0\}}.
\end{align*}
 We then complete the proof of (b) by using Lemma \ref{lem1}.

\end{proof}

\begin{proof}[{\bf Proof of Corollary \ref{th-2}}]
We still consider $\rho_o\in C^2(M)$ for the sake of brevity.
In case (a), for $\sigma\in (0,1)$,
  define the Lyapunov function by
$$W_{\sigma}(r)=\int_{R_0}^{r}h_k(s)^{1-d}\e^{\sigma V(s)}\,\vd s+1,~~\mbox{for all}~~r\geq R_0.$$
Then using a similar calculation as in the proof of Theorem \ref{th1}(a), we have
\begin{align*}
\frac{LW_{\sigma}(\rho_o(x))}{W_{\sigma}(\rho_o(x))}
\leq & -\ff {(1-\sigma)h_k(\rho_o(x))^{1-d}\e^{\sigma V(\rho_o(x))}V'(\rho_o(x))}{W_{\sigma}(\rho_o(x))}
\end{align*}
 for all $\rho_o(x)\geq R_0$. Therefore, there exists a positive constant $b$ such that
\begin{align*}
\frac{LW_{\sigma}(\rho_o(x))}{W_{\sigma}(\rho_o(x))}\leq& -\theta(\rho_o(x))\mathbb{1}_{\{\rho_o(x)\geq R_0\}}+b\mathbb{1}_{\{\rho_o(x)< R_0\}},
\end{align*}
which, together with Lemma \ref{lem1}, implies (a).

In case (b), by assumption \eqref{assumption-b}, we have that for all $\rho_o(x)\geq R_0$,
\begin{align*}
   \delta |\nabla {V}(\rho_o(x))|^2-\Delta {V}(\rho_o(x))&=\delta |{V}'(\rho_o(x))|^2-{V}'(\rho_o(x))\ff{(d-1)h_k'(\rho_o(x))}{h_k(\rho_o(x))}-{V}''(\rho_o(x))> 0.
\end{align*}
 Combining this with Theorem \ref{th-1}(b), we complete the proof of (b).
\end{proof}

\section{Application to convolution  probability measures on $\R^d$}
In this section, we first apply the results in Section 2 to the convolution probability measures
on $\R^d$ and then give the proofs.
\subsection{Main results}
For each $x\in \Rd$, let
$$\nu_x(\vd z)=\frac{1}{p_{\nu}(x)}\e^{-V(x-z)}\,\nu(\vd z).$$
For any non-increasing function $\theta:[0,\infty)\rightarrow (0,\infty)$, let
\begin{equation}\label{beta-phi}
 H_{\theta}(r)= (\mu+\nu)\l(|x|\geq \frac{1}{2}\theta^{-1}\l(1/r\r)\r),\quad r>0.
\end{equation}
By Theorem \ref{th-1}, we have the following first main result.
%For any positive continuous function $\psi$ and positive constant $0<\sigma<1$,
%$$W_{\sigma}^{\psi}(r)=\int_{R_0}^{r}(s^{1-d}\exp[\sigma\int_{R_0}^s\psi(u)\,\vd u])\,\vd s$$
 \begin{theorem}\label{th1}
Let $V\in C^1(\mathbb{R}^d)$ such that $\mu(\vd x)={\e}^{-V(x)}\,\vd x$ is a probability measure on $\Rd$, and let $\nu$ be another probability measure on $\Rd$ such that $p_{\nu}\in C^{1}(\R^d)$.
 \begin{enumerate}
   \item [$(a)$]
   Assume that for some positive constant $R_0$, one has
   \bequ\label{equ-e1}
   \psi(s):=\frac{1}{s}\inf_{|x|=s}\int_{\R^d}\l<\nabla V(x-z), x\r>\nu_x(\vd z)>0,~  ~~s\geq R_0.
   \enqu
  Then for any  $\sigma\in (0,1)$, $\mu\ast \nu$ satisfies the weak Poincar\'{e} inequality with
   $\alpha(r)=cH_{\theta}^{-1}(r)$ for some positive constant $c$, where
   \bequ\label{phi-add}
   \theta(s)=\inf\l\{\ff {(1-\sigma)\psi(r)r^{1-d}\exp[\sigma\int_{R_0}^{r}\psi(u)\,\vd u]}{ \int_{R_0}^{r}t^{1-d}\exp[\sigma\int_{R_0}^t\psi(u)\,\vd u]\,\vd t+1}: r\in [R_0,s\vee R_0]\r\}.
   \enqu
   \item [$(b)$] Let $V\in C^2(\Rd)$ such that for some constant $R_0$ and $\delta\in (0,1)$, one has
   \begin{align}\label{add-cd}
       \theta(s)=(1-\delta)\inf_{|x|\in [R_0, s\vee R_0]}\int_{\Rd}\l(\delta|\nabla V(x-z)|^2-\Delta V(x-z)\r)\nu_x(\vd z)>0.
    \end{align}
    Then $\mu\ast \nu$ satisfies the weak Poincar\'{e} inequality with
       $\alpha(r)=cH_{\theta}^{-1}(r)$ for some positive constant $c$.
 \end{enumerate}
\end{theorem}

\begin{remark}\label{Rem-3}
\begin{itemize}
  \item [(i)] In \cite{WW}, the authors prove that if the function $\phi$ in the Lyapunov condition satisfies
\begin{align}\label{phi-cond}
\liminf_{|x|\rightarrow \infty}\phi(x)=\infty,
\end{align}
then there exists a super Poincar\'{e} inequality with respect to  $\mu*\nu$. Let $\tilde{\theta}(r)=\inf_{|x|\geq r\vee R_0}\phi(x)$. Then \eqref{phi-cond} holds if and only if $\lim_{r\ra \infty} \tilde{\th}(r)=\infty$. Note that in this case, to keep as much information about $\phi$ as possible, it is better for us to choose $\tilde{\th}(|x|)$ instead of $\theta(|x|):=\inf_{|y|\in [R_0,R_0\vee |x|]}\phi(y)$ used in Theorem \ref{th1} to control  $\phi(x)$.
 However, in this article, we take more consideration of the following case for weak Poincar\'{e} inequalities:
  $$\liminf_{|x|\rightarrow \infty}\phi(x)=0.$$
  So in the case \eqref{phi-cond}, we should refer the reader to \cite{WW} for super Poincar\'{e} inequalities.
  \item [(ii)] In the proof of Theorem \ref{th1} (a), it provides a new and reasonable Lyapunov function such that \cite[Theorem 4.1 (a)]{WW} can be improved as follows.
 Recall that $\mu$ satisfies the super Poincar\'{e} inequality with $\beta: (0,\infty)\rightarrow (0,\infty)$ if
 \begin{align*}
 \mu(f^2)\leq r\mu(|\nabla f|^2)+\beta(r)\mu(|f|)^2,\quad r>0,\ f\in C_b^1(\R^d).
 \end{align*}

\noindent {\bf Theorem A.} {\it
 Let $V\in C^1(\mathbb{R}^d)$ such that $\mu(\vd x)={\e}^{-V(x)}\,\vd x$ is a probability measure on $\Rd$, and let $\nu$ be another probability measure on $\Rd$ such that $p_{\nu}\in C^{1}(\R^d)$.
Let $$\xi(r,s)=\l(1+s^{-\frac{d}{2}}\r)\frac{\sup_{|x|\leq r}\e^{-(\frac{d}{2}+1)V(x)}}{\inf_{|x|\leq r}\e^{-(\frac{d}{2}+2)V(x)}}.$$
If $\liminf_{r\rightarrow \infty}\psi(r)=+\infty$, where $\psi$ is defined as in \eqref{equ-e1},
then $\mu*\nu$ satisfies the super Poincar\'{e} inequality with
$$\beta(r)=c(1+\xi\l(\gamma\l({2}/{r}\r), {r}/{2})\r),$$
where $c$ is  some positive constant and
$$\gamma(s):=\inf\l\{t>0:\inf_{r\geq t\vee R_0}\ff {(1-\sigma)\psi(r)r^{1-d}\exp[\sigma\int_{R_0}^{r}\psi(u)\,\vd u]}{ \int_{R_0}^{r}s^{1-d}\exp[\sigma\int_{R_0}^s\psi(u)\,\vd u]\,\vd s+1}\geq s\r\}$$
for any $\sigma\in (0,1)$ and  some positive constant $R_0$ such that $\psi(r)>0$ for all $r\geq R_0$.
}

In the following subsection, we will give a brief explanation for the proof of this theorem and use the example in \cite[Theorem 4.4]{WW} to show the benefit of this result.
\end{itemize}
\end{remark}

From \eqref{equ-e1} and \eqref{add-cd}, it is easy to see that if the function $p_{\nu}$ has previous estimates, then Theorem \ref{th1} can be simplified as follows.

 \begin{theorem}\label{th2}
Let $\mu(\vd x)={\e}^{-V(x)}\,\vd x$ be a probability measure on $\Rd$ and  $\nu$ be another probability measure on $\Rd$ such that $p_{\nu}\in C^1(\R^d)$. Set
$$\e^{-\tilde{V}_{\nu}(s)}=\inf_{|x|=s}\int_{\R^d} \e^{-V(x-z)}\,\nu(\vd z)\ \ \mbox{and}\ \ \e^{-\hat{V}_{\nu}(s)}=\sup_{|x|=s}\int_{\R^d} \e^{-V(x-z)}\,\nu(\vd z).$$
 \begin{enumerate}
   \item [$(a)$] If $\tilde{V}_{\nu}\in C^{1}([0,\infty))$ such that for some positive constant $R_0$ and any $\sigma\in (0,1)$, one has
   \begin{align}\label{p-sigma}
     \theta(s):=\inf_{r\in[R_0, s\vee R_0]}\ff{(1-\sigma)\tilde{V}_{\nu}'(r){\e}^{\sigma \tilde{V}_{\nu}(r)}r^{1-d}} {\int_{R_0}^{r}s^{1-d}{\e}^{\sigma \tilde{V}_{\nu}(s)} \,\vd s +1}>0,
   \end{align}
then $\mu*\nu$ satisfies the  weak Poincar\'{e} inequality with
\begin{equation}\label{beta-phi-2}\alpha(r)=c\inf\l\{\frac{\sup_{0\leq t< 2s}\e^{\tilde{V}_{\nu}(t)-\hat V_{\nu}(t)}}{\theta(2s)}: (\mu+\nu)\l(|x|\geq s\r)\leq r,\ s>0\r\}\end{equation}
 for some positive constant $c$.

   \item [$(b)$] If $\tilde{V}_{\nu}\in C^2([0, \infty))$ such that for some positive constant $R_0$ and $\delta \in (0,1)$,
   \begin{align}\label{add-cd-1}
      \theta(s):=(1-\delta)\inf_{r\in[R_0, s\vee R_0]}\l[\delta |\tilde{V}_{\nu}'(r)|^2-\tilde{V}_{\nu}'(r)\frac{d-1}{r}-\tilde{V}_{\nu}''(r)\r]>0,
    \end{align}
        then $\mu\ast \nu$ satisfies the weak Poincar\'{e} inequality with
       $\alpha(r)$ defined as in \eqref{beta-phi-2} for some positive constant $c$.
 \end{enumerate}
\end{theorem}

Next, we shall apply above results to the convolution with compactly supported probability measures. Note that, if $\nu$ is a probability measure with compact support, then the function  $p_{\nu}$ is obviously differentiable on $\R^d$. Thus, by Theorem \ref{th1}, we  obtain the following corollary directly.
 \begin{corollary}\label{cor1}
Let $V\in C^1(\mathbb{R}^d)$ such that $\mu(\vd x)={\e}^{-V(x)}\,\vd x$ is a probability measure on $\Rd$ and let $\nu$ be another probability measure on $\Rd$ with $R:=\sup\{|z|:z\in {\rm supp} \nu\}<\infty$.
\begin{enumerate}
  \item [$(a)$]  Assume that for some  positive  constant $R_0>R$, one has
    \begin{equation}\label{eq-111}
      \psi(s):=\frac{1}{s}\inf_{s-R\leq |u|\leq R+s}(\l<u,\nabla V(u)\r>-R|\nabla V(u)|)>0,~~~s\geq R_0.
    \end{equation}

 For any  $\sigma\in (0,1)$, let
   \bequ\label{phi}
   \theta(s)=\inf\l\{\ff {(1-\sigma)\psi(r)r^{1-d}\exp[\sigma\int_{R_0}^{r}\psi(u)\,\vd u]}{ \int_{R_0}^{r}(s^{1-d}\exp[\sigma\int_{R_0}^s\psi(u)\,\vd u])\,\vd s+1}: r\in [R_0,s\vee R_0]\r\}.
   \enqu
  Then $\mu\ast \nu$ satisfies the weak Poincar\'{e} inequality with
   $\alpha(r)=cH_{\theta}^{-1}(r)$ for some positive constant $c$.

  \item [$(b)$]Let $V\in C^2(\Rd)$ such that for some positive constants $R_0>R$ and $\delta \in (0,1)$, one has
      \begin{align}\label{theta-2}
       \theta(s):=(1-\delta)\inf_{R_0-R\leq |u|\leq R+s\vee R_0}(\delta|\nabla V(u)|^2-\Delta V(u))>0.
      \end{align}
  Then $\mu\ast\nu$ satisfies the weak Poincar\'{e} inequality with $\alpha(r)=cH_{\theta}^{-1}(r)$
  for some positive constant $c$.

\end{enumerate}
\end{corollary}
\begin{remark}\label{rem-1}
We remark that  in Corollary \ref{cor1}, due to the compactness of $\nu$ and the monotonicity of $\theta$, there exists a positive constant $r_0$ such that
$$H_{\theta}(r)=\mu\l(|x|\geq \frac{1}{2}\theta^{-1}\big(1/r\big)\r), ~~~r\in (r_0,\infty).$$
\end{remark}

 By Theorem \ref{th2}, we  obtain the following corollary directly.
 \begin{corollary}\label{cor2}
Let $V\in C^1(\mathbb{R}^d)$ such that $\mu(\vd x)={\e}^{-V(x)}\,\vd x$ is a probability measure on $\Rd$ and let $\nu$ be another probability measure on $\Rd$ with $R:=\sup\{|z|:z\in {\rm supp} \nu\}<\infty$. For some constant  $R_0>R$ and any $s\geq R_0$, let
$$\tilde{V}(s)=\sup_{s-R\leq |x|\leq R+s}V(x)\ \ \mbox{and}\ \ \hat{V}(s)=\inf_{s-R\leq |x|\leq R+s}V(x). $$
Then the assertions in Theorem \ref{th2}(a)(b) still hold by replacing $\tilde{V}_{\nu}$ and $\hat{V}_{\nu}$ with $\tilde{V}$ and $\hat{V}$, respectively.
\end{corollary}

\subsection{Proofs}

Using Lemma \ref{lem1}, we  complete the proof of Theorem \ref{th1}.

\noindent{\bf Proof of Theorem \ref{th1}.} Let $L_\nu=\De-\nn V_\nu$. First, if the Lyapunov condition {\bf(L)} holds for $L_\nu$ with some function $\phi=\theta(|\cdot|)$, where $\theta:(0,\infty)\rightarrow [0,\infty)$ is a non-increasing function, then we have
\begin{align*}
&\mu\ast\nu\l(\theta(|x|)\leq \frac{1}{r}\r)\\
&\leq  \mu\ast\nu(|x|\geq \theta^{-1}(1/r))
  =\iint_{\{|x|\geq \theta^{-1}(1/r)\}}\e^{-V(x-z)}\,\nu(\vd z)\,\vd x\\
&\leq \iint_{\{|x-z|\geq \frac{1}{2}\theta^{-1}(1/r)\}}\e^{-V(x-z)}\,\nu(\vd z)\,\vd x+\iint_{\{|z|\geq \frac{1}{2}\theta^{-1}(1/r)\}}\e^{-V(x-z)}\,\nu(\vd z)\,\vd x\\
&= \iint_{\{|x-z|\geq \frac{1}{2}\theta^{-1}(1/r)\}}\e^{-V(x-z)}\,\vd x\,\nu(\vd z)+\int_{\{|z|\geq \frac{1}{2}\theta^{-1}(1/r)\}}\int\e^{-V(x-z)}\,\vd x\,\nu(\vd z)\\
&=\mu\l(2|x|\geq \theta^{-1}(1/r)\r)+\nu\l(2|z|\geq \theta^{-1}(1/r)\r)\\
&=H_{\theta}(r).
\end{align*}
Hence, by Lemma \ref{lem1}, $\mu*\nu$ satisfies a weak Poincar\'{e} inequality with $\alpha(r)=cH_{\theta}^{-1}(r)$ for some positive constant $c$. We now turn to construct some suitable Lyapunov functions.

In case (a),
  define
$$W_{\sigma}(r)=\int_{R_0}^{r}\Big(s^{1-d}\exp\Big[\sigma\int_{R_0}^s\psi(u)\,\vd u\Big]\Big)\,\vd s+1,~~\mbox{for all}~~r\geq R_0,$$
where $0<\sigma<1$. Then by a similar discussion as in the proof of Corollary \ref{th-2} for $k=0$ and $L_{\nu}$,  we have
that  there exists a constant $b>0$ such that
\begin{align}
\frac{L_{\nu}W_{\sigma}(|x|)}{W_{\sigma}(|x|)}\leq&
 -\ff {(1-\sigma)\psi(|x|)|x|^{1-d}\exp\big[\sigma\int_{R_0}^{|x|}\psi(u)\,\vd u\big]}{ \int_{R_0}^{|x|}s^{1-d}\exp\big[\sigma\int_{R_0}^s\psi(u)\,\vd u\big]\,\vd s+1}\mathbb{1}_{\{|x|\geq R_0\}}+b\mathbb{1}_{\{|x|< R_0\}}\label{phi-new}\\
\leq & -\theta(|x|)\mathbb{1}_{\{|x|\geq R_0\}}+b\mathbb{1}_{\{|x|< R_0\}}.\notag
\end{align}

In case (b), we consider a smooth function such that $W(x)=\e^{(1-\delta)V_\nu(x)}$ for $|x|\geq R_0$ and $W(x)\geq 1$ for all $x\in \Rd$. Using the same argument as in the proof of Theorem \ref{th-1}, we have
\begin{align}\label{eq-2}
\frac{L_{\nu}W(x)}{W(x)}&\leq-(1-\delta)(\delta |\nabla V_{\nu}|^2-\Delta V_{\nu})\mathbb{1}_{\{|x|\geq R_0\}}+b\mathbb{1}_{\{|x|< R_0\}}.
\end{align}
Moreover,  for any $|x|\geq R_0$,
\begin{align}\label{eq-3}
    \delta |\nabla V_{\nu}(x)|^2-\Delta V_{\nu}(x)&=\int_{\Rd}(|\nabla V(x-z)|^2-\Delta V(x-z))\nu_{x}(\vd z)-(1-\delta)|\nabla V_{\nu}(x)|^2\notag\\
    &\geq \int_{\Rd}(\delta |\nabla V(x-z)|^2-\Delta V(x-z))\nu_x(\vd z)\notag\\
    &\geq \frac{1}{1-\delta}\theta(|x|).
\end{align}
 Combining this with \eqref{eq-2}, we complete the proof of (b).
\endzm

\begin{proof}[{\bf Proof of Theorem A}]
Let $$\theta(r)=\ff {(1-\sigma)\psi(r)r^{1-d}\exp[\sigma\int_{R_0}^{r}\psi(u)\,\vd u]}{ \int_{R_0}^{r}s^{1-d}\exp[\sigma\int_{R_0}^s\psi(u)\,\vd u]\,\vd s+1}.$$
By \eqref{phi-new}, we know that $L_{\nu}$ satisfies
$$\frac{L_{\nu}W_{\sigma}(|x|)}{W_{\sigma}(|x|)}
\leq  -\theta(|x|)\mathbb{1}_{\{|x|\geq R_0\}}+b\mathbb{1}_{\{|x|< R_0\}}.$$
It is easy to see that
if $\liminf_{r\rightarrow \infty}\psi(r)=+\infty$, then $\liminf_{r\rightarrow \infty}\theta(r)=+\infty$.
Thus by \cite[Lemma 4.2]{WW}, we complete the proof directly.
\end{proof}
Now we use Theorem A to prove the following result.
\begin{example}
Let $V(x)=c+|x|^p$ for some $p>1$ and $c\in \R$ such that $\mu(\vd x):=\e^{-V(x)}\,\vd x$ is a probability measure
on $\R^d$. Let $\nu$ be any compactly supported probability measure. Then there exists a constant $c>0$ such that $\mu*\nu$ satisfies the super Poincar\'{e}
inequality with
$$\beta(r)=\exp(cr^{-\frac{p}{2(p-1)}}),\quad r>0. $$
\end{example}
\begin{proof}[{\bf Proof}]
Suppose that $\nu$ is supported on $\{x:|x|\leq R\}$ for some positive constant $R$. Then
\begin{align}\label{psi-est}
\psi(s)\geq \frac{1}{s}\inf_{|x|=s, |z|\leq R}\l<\nabla V(x-z), x\r>\geq \frac{p|s-R|^{p-1}s}{s+R}.
\end{align}
Thus there exists a positive constant $R_0$, for $r>R_0$,
$$\theta(r)=\ff {(1-\sigma)\psi(r)r^{1-d}\exp\big[\sigma\int_{R_0}^{r}\psi(u)\,\vd u\big]}{ \int_{R_0}^{r}s^{1-d}\exp\big[\sigma\int_{R_0}^s\psi(u)\,\vd u\big]\,\vd s+1}\geq c_1r^{2(p-1)}$$
for some positive constant $c_1$. Thus $\gamma(u)\leq c_2(1+u^{\frac{1}{2(p-1)}}), u>0$
holds for some positive constant $c_2$. Moreover, as explained in the proof of \cite[Example 4.4]{WW}, one has
$$\xi(r,s)\leq c_3(1+s^{-d/2})\e^{c_4r^p},\quad s,r>0$$
 for some positive constants $c_3,c_4$. So the desired assertion follows by using Theorem A.
\end{proof}
However, by \cite[Theorem 4.1 (a)]{WW}, it is easy to calculate that $\mu*\nu$ satisfies
the super Poincar\'{e} inequality with
$$\beta(r)=\exp(cr^{-\frac{p}{p-1}}),\quad r>0,$$
which is less sharp than that presented in this example.

Let us continue with the proofs of main results in Subsection 3.1.

\noindent{\bf Proof of Theorem \ref{th2}.} Let $\tilde{L}=\De-\nn \tilde{V}_{\nu}$ and $\tilde{\mu}(\vd x)=\e^{-\tilde{V}_{\nu}(|x|)}\,\vd x$. First, if the Lyapunov condition {\bf(L)} holds for $\tilde{L}$ with some function $\phi(\cdot)=\theta(|\cdot|)$, where $\theta$ is a positive and non-decreasing function on $\R^+$, then for any $f\in C_b^1(\mathbb{R}^d)$ with $\mu(f)=0$ and $x_0\in \mathbb{R}^d$ such that $f(x_0)=\frac{1}{\tilde{\mu}(B_{r_0})}(\int_{B_{r_0}}f\,\vd \tilde{\mu})$,
\begin{align*}
\mu*\nu(f^2)&\leq \inf_{c\in \mathbb{R}}\mu*\nu(f-c)^2\\
&\leq \int_{\phi>1/s}(f-f(x_0))^2\,\vd \mu*\nu +\int_{\phi\leq 1/s}(f-f(x_0))^2\,\vd \mu*\nu\\
&\leq s\int_{\phi>1/s}\phi (f-f(x_0))^2\,\vd \mu*\nu +\mu*\nu \l(\phi \leq 1/s\r){\rm Osc}(f)^2\\
&\leq s\sup_{\phi>1/s}\e^{\tilde{V}_{\nu}(|x|)-\hat V_{\nu}(|x|)}\int \phi (f-f(x_0))^2\,\vd \tilde{\mu}+\mu*\nu\l(\phi\leq 1/s\r){\rm Osc}(f)^2\\
&\leq s\sup_{0\leq t< \theta^{-1}\l(1/s\r)}\e^{\tilde{V}_{\nu}(t)-\hat V_{\nu}(t)}\int \phi (f-f(x_0))^2\,\vd \mu*\nu\\
&\quad+(\mu+\nu)\l(|x|\geq \frac{1}{2}\theta^{-1}\l(1/s\r)\r){\rm Osc}(f)^2,\qquad \quad \ 1/s> \inf \theta.
\end{align*}
Let $r=\frac{1}{2}\theta^{-1}\l(1/s\r)$. Then, using a similar argument as in the inequality \eqref{estimate-11} we obtain
\begin{align*}
\mu*\nu(f^2)
&\leq \frac{\sup_{0\leq t< 2r}\e^{\tilde{V}_{\nu}(t)-\hat V_{\nu}(t)}}{\theta(2r)}\int |\nabla f|^2\,\vd \mu*\nu+(\mu+\nu)\l(|x|\geq r\r){\rm Osc}(f)^2,\qquad \quad \ r>0.
\end{align*}
It follows that $\mu*\nu$ satisfies a  weak Poincar\'{e} inequality with
$$\alpha(s)=c\inf\l\{\frac{\sup_{0\leq t< 2r}\e^{\tilde{V}_{\nu}(t)-\hat V_{\nu}(t)}}{\theta(2r)}: (\mu+\nu)\l(|x|\geq r\r)\leq s,\ r>0\r\}$$
 for some positive constant $c$.
Now it suffices for us to construct some suitable Lyapunov functions.

In case (a),
  define the Lyapunov function by
$$W_{\sigma}(|x|)=\int_{R_0}^{|x|}s^{1-d}\e^{\sigma \tilde{V}_{\nu}(s)}\,\vd s+1,~~\mbox{for all}~~|x|\geq R_0.$$
Then by a similar calculation as in the proof of Corollary \ref{th-2}\,(a), we have that for all $|x|\geq R_0$,
 there exists a positive constant $b$ such that
\begin{align*}
\frac{\tilde{L}W_{\sigma}(|x|)}{W_{\sigma}(|x|)}\leq& -\theta(|x|)\mathbb{1}_{\{|x|\geq R_0\}}+b\mathbb{1}_{\{|x|< R_0\}}.
\end{align*}

In case (b), we consider a smooth function such that $W(x)=c\e^{(1-\delta)\tilde{V}_{\nu}(|x|)}$ for $|x|\geq R_0$. Then,
\begin{align*}
\frac{\tilde{L}W(x)}{W(x)}&\leq-(1-\delta)\l[\delta |\tilde{V}'_{\nu}(|x|)|^2-\tilde{V}_{\nu}'(|x|)\ff{d-1}{|x|}-\tilde{V}''_{\nu}(|x|)\r]\mathbb{1}_{\{|x|\geq R_0\}}+b\mathbb{1}_{\{|x|< R_0\}}\\
&=-\theta(|x|)\mathbb{1}_{\{|x|\geq R_0\}}+b\mathbb{1}_{\{|x|< R_0\}}.\ \ \hspace{7cm} \Box
\end{align*}

\medskip

\noindent{\bf Proof of Corollary \ref{cor1}.}
In case (a). It is easy to see that
\begin{ews*}
&\int_{\Rd}\l<x,\nabla V(x-z)\r>\nu_x(\vd z)\\
&\quad=\int_{\Rd}\big(\l<x-z,\nabla V(x-z)\r>+\l<z,\nabla V(x-z)\r>\big)\nu_x(\vd z)\\
&\quad\geq \int_{\Rd}\big(\l<x-z,\nabla V(x-z)\r>-R|\nabla V(x-z)|\big)\nu_x(\vd z)\\
&\quad= \int_{\{|z|\leq R\}} \big(\l<x-z,\nabla V(x-z)\r>-R|\nabla V(x-z)|\big)\nu_x(\vd z).
\end{ews*}
Then according to the definitions of $\psi$, we have that for any $s\geq R_0(>R)$,
\begin{ews*}
&\inf_{|x|=s}{\int_{\Rd}\l<x,\nabla V(x-z)\r>\nu_x(\vd z)} \\
&\quad\geq \inf_{|x|=s}\int_{\{|z|\leq R\}}\l({\l<x-z,\nabla V(x-z)\r>-R|\nabla V(x-z)|}\r) \nu_x(\vd z)\\
&\quad\geq  \inf_{s-R\leq |u|\leq s+R}(\l<u,\nabla V(u)\r>-R|\nabla V(u)|)\\
&\quad = s\psi(s).
\end{ews*}
Then, we  complete the proof of (a) due to Theorem \ref{th1}(a).

In case (b). For any $s\geq R_0\,(>R)$, we have that for $s\geq R_0$,
\begin{ews*}
&\inf_{R_0\leq |x|\leq s}\int_{\Rd}\Big(\delta|\nabla V|^2(x-z)-\Delta V(x-z)\Big)\nu_{x}(\vd z)\\
&\quad=\inf_{R_0\leq |x|\leq s}\int_{\{|z|\leq R\}}\Big(\delta|\nabla V|^2(x-z)-\Delta V(x-z)\Big)\nu_{x}(\vd z)\\
&\quad\geq\inf_{(R_0-R)\leq |u|\leq R+s}\l(\delta|\nabla V|^2(u)-\Delta V(u)\r)>0,
\end{ews*}
which leads to complete the proof by Theorem \ref{th1}(b).
\endzm

\noindent{\bf Proof of Corollary \ref{cor2}.}
The results follow from Theorem \ref{th2} and the following fact: there exists constant $R_0>R$ such that for $s\geq R_0$,
\begin{align*}
&\sup_{|x|=s}\int_{\R^d}\e^{-V(x-z)}\,\nu(\vd z)\leq \sup_{|x|=s}\sup_{|z|\leq R}\e^{-V(x-z)}\leq  \sup_{s-R\leq|u|\leq s+R} \e^{-V(u)}=\e^{-\tilde{V}(s)}; \\
&\inf_{|x|=s}\int_{\R^d}\e^{-V(x-z)}\,\nu(\vd z)\geq \inf_{|x|=s}\inf_{|z|\leq R}\e^{-V(x-z)}\geq  \inf_{s-R\leq|u|\leq s+R} \e^{-V(u)}=\e^{-\hat{V}(s)}.\ \ \Box
\end{align*}

\section{Examples}
\hspace{0.4cm}
In this section, we present the following examples to illustrate the results obtained in Section 3.
As an application of  Theorem \ref{th2}, we present  below an example where the support of $\nu$ is unbounded and disconnected.
\beg{example}\label{ex-0}
Let $d=1$. For   $0<\delta<1$ and $p>0$, let  $V(x)=c+(1+x^2)^{\ff \delta 2}$ and
 $$\nu(\d z)=\ff 1 \ga \sum_{i\in \Z}\ff {\de_i(\d z)} {1+|z|^{1+p}},$$
  where
$$c=\log \int_{\R}\e ^{-(1+x^2)^{\ff \delta 2}}\,\d x~~\mbox{and}~~\ga=\sum_{i\in \Z}\ff 1 {1+|i|^{1+p}}.$$
Then  there exists a positive constant $C$ such that the weak Poincar\'e inequality for $\mu*\nu$ holds with $\al(s)=Cs^{-2/p}$ for all $s>0$.
\end{example}

\zm
We use Theorem 3.3(a) to give the proof. First we need to estimate $p_{\nu}$.
It is easy to see that
$$p_{\nu}(x)=\frac{\e^{-c}}{\gamma}\sum_{i\in \Z}\ff {\e^{-[1+(x-i)^2]^{\ff {\de} {2}}}}{1+|i|^{1+p}}=\frac{\e^{-c}}{\gamma}\sum_{k\in \Z}\ff {\e^{-[1+((x)-k)^2]^{\ff {\de} {2}}}}{1+|[x]-k|^{1+p}},$$
where $x=[x]+(x)$ and $[x]$ is the integral part of $x$.
Moreover, as
$$\frac{1}{2}k^2\leq 1+((x)-k)^2\leq 2+k^2,$$
we have
\begin{align}\label{com-1}
\sum_{k\in \Z}\ff {\e^{-(2+k^2)^{\ff {\de} {2}}}}{1+|[x]-k|^{1+p}}\leq \sum_{k\in \Z}\ff {\e^{-[1+((x)-k)^2]^{\ff {\de} {2}}}}{1+|[x]-k|^{1+p}}\leq \sum_{k\in \Z}\ff {\e^{-(\frac{1}{2})^{\ff{\de}{2}}k^\de}}{1+|[x]-k|^{1+p}}.
\end{align}
To deal with the terms on the both sides of the inequality above, we need  the following  estimates:
$$\ff {|[x]|^{p+1}} {1+|[x]-k|^{p+1}}\leq \ff {2^p\Big(|[x]-k|^{p+1}+|k|^{p+1}\Big)} {1+|[x]-k|^{p+1}}\leq 2^p\Big(1+|k|^{p+1}\Big).$$
Using these inequalities and the dominated convergence theorem, we have
$$\lim_{|x|\ra+\infty}\sum_{k\in\Z}\ff {|[x]|^{p+1}\e^{-(\frac{1}{2})^{\ff{\de}{2}}k^{\de}}} {1+|[x]-k|^{p+1}}=\sum_{k\in\Z}\e^{-(\frac{1}{2})^{\ff{\de}{2}}k^{\de}};$$
and
$$\lim_{|x|\ra+\infty}\sum_{k\in\Z}\ff {|[x]|^{p+1}\e^{-(2+k^2)^{\ff {\de} {2}}}} {1+|[x]-k|^{p+1}}=\sum_{k\in\Z}\e^{-(2+k^2)^{\ff {\de} {2}}}.$$
 Combining these with \eqref{com-1} yields
 $$p_{\nu}(x)\asymp \frac{1}{|x|^{1+p}}.$$
Here and in what follows, for any functions $f$ and $g$, we write ``$f \asymp g$" if there exist positive constants $c_1$ and $c_2$ such that $c_2f\leq g\leq c_1f$.
 It then follows that
 $$\hat{V}_{\nu}(s)=\tilde{V}_{\nu}(s)=\log(1+s)^{1+p}+{\rm o}(\log(1+s)).$$
  By this and the definition of $\theta$ in \eqref{p-sigma},
  there exists some positive constant $R_0$  such that
  \begin{align}\label{theta-1}
     \theta(s):=\inf_{r\in[R_0, s\vee R_0]}\ff{(1-\sigma)\tilde{V}_{\nu}'(r){\e}^{\sigma \tilde{V}_{\nu}(r)}r^{1-d}} {\int_{R_0}^{r}t^{1-d}{\e}^{\sigma \tilde{V}_{\nu}(t)} \,\vd t +1}\asymp \frac{1}{s^2}.
   \end{align}
   Moreover, it is easy to calculate that for large $r$,
   $$(\mu+\nu)(|x|\geq r)\asymp r^{-p}. $$
   By this and \eqref{theta-1}, we conclude that there exists some positive constant $C$ such that
$$H_\theta(r)\leq Cr^{p/2},$$
 which completes the proof by  Theorem \ref{th2}\,(a).
\endzm

\begin{example}\label{em1}
Let $V(x)=c+|x|^{p}$ for some $0<p<1$, and $\mu(\vd x)=\e^{-V(x)}\,\vd x$. Then for any probability measure $\nu$ with $R:=\sup\{|z|: z\in {\rm supp}\ \nu\}<\infty$,  there exists some positive constant $C$ such that the weak poincar\'{e} inequality for $\mu*\nu$ holds with
$$\alpha(s)=C\Big[1+\log\l(1+\ff 1 s\r)\Big]^{\ff {2(1-p)} {p}},\quad s>0.$$
\end{example}
\zm a) \emph{Method 1.}\
  It is easy to see that
 $$\inf_{|x|=s}\l(\l<\nabla V(x),x\r>-R|\nabla V(x)|\r)=ps^{p-1}(s-R).$$
Then there exists $R_0>R$ such that for $|x|\geq R_0$,
\begin{align*}
\inf_{|x|-R\leq s\leq |x|+R}ps^{p-1}(s-R)\asymp |x|^{p}.
\end{align*}
Thus, we can choose  $\psi(|x|):= c|x|^{p-1}$ and then have that for $|x|\geq R_0$,
\beg{ews*}
\ff{c(1-\sigma)|x|^{p-1}|x|^{1-d}\e^{c\si |x|^{p}}} {\int_{R_0}^{|x|} u^{1-d}\e^{c\si u^{p}}\,\d u+1} \geq C\ff {|x|^{p-d}\e^{c\si |x|^{p}}} {|x|^{2-d-p}\e^{c\sigma |x|^{p}}}=C|x|^{2(p-1)}.
\end{ews*}
It follows from the definition of $\theta$ in \eqref{phi} that
 $$\theta(|x|)\asymp  |x|^{2(p -1)}, ~~ \mbox{for all}~~ |x|\geq R_0. $$
From this, we obtain that for any $r>0$,
\beg{ews*}
H_{\theta}(r)&=\mu\Big(2|x|\geq \theta^{-1}(1/r)\vee R_0\Big)\leq C\int_{cr^{\ff 1 {2(1-p)}}}^{\infty}\e^{-u^p}u^{d-1}\,\d u\leq C\e^{-cr^{\ff {p} {2(1-p)}}}r^{\ff {d-p} {2(1-p)}}.
\end{ews*}
Now using Corollary \ref{cor1} (a),  we conclude that there exists some positive constant $C$   such that
$$\al(s)=C\Big[1+\log\l(1+\ff 1 s\r)\Big]^{\ff {2(1-p)} {p}}.$$

b) \emph{ Method 2.}\  It is easy to calculate that for $\delta>0$ and $|x|>0$,
\begin{align*}
{\delta|\nabla V(x)|^2-\Delta V(x)=\delta p^2|x|^{2(p-1)}-p(d+p-2)|x|^{p-2}.}
\end{align*}
Thus, there exists some constant $R_0>0$ such that  for all $|x|\geq R_0$,
$$\inf_{ |u|\leq |x|+R}|u|^{2(p-1)}\geq  (|x|+R)^{2(p-1)}.$$
So the function $\theta$ in Corollary \ref{cor1} (b) satisfies
$$\theta(r)\asymp r^{2(p-1)},~r\geq R_0.$$
The rest of the proof is similar by using Corollary \ref{cor1} (b), so we omit it.
\endzm

Next, the following examples are to illustrate Corollary \ref{cor2}.
\begin{example}\label{add-ex-2}
For $p>0$, let $V(x)=c+(d+p)\log (1+|x|)$. Then for any  probability measure {$\nu$ with $R:=\sup\{|z|: z\in {\rm supp}\ \nu\}<\infty$}, there exists some positive constant $C$ such that $\mu*\nu$ satisfies the weak poincar\'{e} inequality  with
$$\al(s)=Cs^{-\ff 2p},\quad s>0.$$
\end{example}
\zm We  use Corollary \ref{cor2} to give the proof.
It is easy to see that for $s>R$,
\begin{align*}
  &\tilde{V}(s)=\sup_{s-R\leq |x|\leq R+s}V(x)=c+(d+p)\log(1+R+s),\\
  &\hat{V}(s)=\inf_{s-R\leq |x|\leq R+s}V(x)=c+(d+p)\log(1+s-R).
\end{align*}
Thus, there exists a positive constant $C$ such that
$$\e^{\tilde{V}(s)-\hat{V}(s)}=\frac{(1+R+s)^{d+p}}{(1+s-R)^{d+p}}\leq C,\quad s>R.$$
Moreover, for $\si\in(\ff {d-2} {d+p}\vee 0,1)$, let $\theta$ be in \eqref{phi}. Then there exists a positive constant $R_0>R$ such that
$$\theta(r)=\ff{c(1-\sigma)(d+p)(1+R+r)^{\sigma(d+p)-1}r^{1-d}}{\int _{R_0}^{r}(1+R+s)^{\sigma(d+p)}s^{1-d}\,\vd s+1}\leq c(1+r)^{-2}, ~~~ r\geq R_0\,(>R).$$
Therefore, by Corollary \ref{cor2}, we obtain the results directly. This result also can be proved in a similar way by using $\theta$ constructed in \eqref{theta-2} and Corollary \ref{cor2}.
\endzm

Similarly, we have

\begin{example}\label{add-ex-3}
Let $p>1$ and $V(x)=c+d\log(1+|x|)+p\log\log(\e+|x|)$. Then for any probability measure $\nu$ with $R:=\sup\{|z|: z\in {\rm supp}\ \nu\}<\infty$, there exist  some positive constants $c_1,c_2$  such that the weak poincar\'{e} inequality holds for $\mu*\nu$  with
$$\alpha(r)=c_1\exp[c_2r^{-1/(p-1)}],\quad r>0.$$
\end{example}

\begin{remark}\label{Rem2}
When $\nu=\delta_0$, i.e. $R=0$, Examples \ref{add-ex-2}--\ref{add-ex-3} have been treated in \cite{Wbook1}. Compared with the results in \cite{Wbook1},   the results presented above are more precise. We would like to indicate that  by \cite[Corollary 4.2.2 (1)]{Wbook1}, the $\alpha$ in Example \ref{add-ex-2} implies
the exact main order of $\mu*\nu(|x|>N)$ as $N\rightarrow \infty$. Hence,  using Lyapunov conditions seems to be able to get  better convergence or decay rates for diffusion processes.
\end{remark}

\vspace{0.5cm}
\noindent\textbf{Acknowledgements}  \  The authors would like to thank Professor Feng-Yu Wang for his guidance. The first author was supported by Fonds National
de la Recherche Luxembourg (Open project O14/7628746 GEOMREV), NSFC (Grant No.~A011002)
and Zhejiang Provincial Natural Science Foundation of China (Grant
No. LQ16A010009).

\bigskip

\footnotesize{
\noindent Li-Juan Cheng

\medskip

\noindent Department of Applied Mathematics, Zhejiang University of Technology, Hangzhou 310023,   China

\smallskip

\noindent Mathematics Research Unit, FSTC, University of Luxembourg, Luxembourg, Grand Duchy of
    Luxembourg

\noindent{\it E-mail:} \texttt{chenglj@mail.bnu.edu.cn}, \texttt{lijuan.cheng@uni.lu}

\bigskip

\noindent Shao-Qin Zhang (Corresponding author)

\medskip

\noindent School of Statistics and Mathematics, Central University of Finance and Economics, Beijing, 100081, China

\smallskip

\noindent{\it E-mail:} \texttt{zhangsq@cufe.edu.cn}
}

\end{document}